\documentclass[12pt]{amsart}
\usepackage{amsmath,amssymb,amsbsy,amsfonts,amsthm,latexsym,
                        amsopn,amstext,amsxtra,euscript,amscd,mathrsfs,color,bm}
                        
\usepackage{float} 
\usepackage[english]{babel}
\usepackage{mathtools}

\usepackage[backref=page]{hyperref}
\hypersetup{nesting=true,debug=true,naturalnames=true}
\usepackage{upref}
\usepackage{nicefrac}

\newtheorem{theorem}{Theorem}

\newtheorem{corollary}[theorem]{Corollary}

\newtheorem{conj}[theorem]{Conjecture}

\newtheorem{remark}[theorem]{Remark}

\numberwithin{equation}{section}
\numberwithin{theorem}{section}
\numberwithin{table}{section}
\numberwithin{figure}{section}

\allowdisplaybreaks

\newenvironment{notation}[0]{%
  \begin{list}%
    {}%
    {\setlength{\itemindent}{0pt}
     \setlength{\labelwidth}{1\parindent}
     \setlength{\labelsep}{\parindent}
     \setlength{\leftmargin}{2\parindent}
     \setlength{\itemsep}{0pt}
     }%
   }%
  {\end{list}}

  {\end{list}}

\renewcommand{\a}{\alpha}
\renewcommand{\b}{\beta}
\newcommand{\g}{\gamma}

\newcommand{\Scal}{{\mathcal S}}
\newcommand{\Tcal}{{\mathcal T}}
\newcommand{\Ucal}{{\mathcal U}}
\newcommand{\Vcal}{{\mathcal V}}

\newcommand{\KK}{\mathbb{K}}
\newcommand{\LL}{\mathbb{L}}

\newcommand{\PP}{\mathbb{P}}
\newcommand{\QQ}{\mathbb{Q}}

\newcommand{\ZZ}{\mathbb{Z}}

\newcommand{\Div}{\operatorname{Div}}

\newcommand{\hhat}{{\hat h}}

\newcommand{\logplus}{\log^{\scriptscriptstyle+}}

\newcommand{\ord}{\operatorname{ord}}
\newcommand{\Orbit}{\mathcal{O}}
\newcommand{\Per}{\operatorname{Per}}

\newcommand{\PrePer}{\operatorname{PrePer}}

\renewcommand{\setminus}{\smallsetminus}

\newcommand{\Support}{\operatorname{Supp}}

\newcommand{\Wander}{\operatorname{Wander}}

\newcommand{\ov}[1]{\overline{#1}}
\newcommand{\RSf}{R_{\Scal_{f}}}
\newcommand{\RTf}{R_{\Tcal_{f}}}

\newcommand{\Case}[2]{\smallskip\paragraph{\textbf{\boldmath Case {#1}: {#2}.}}\hfil\break\ignorespaces}

\newcommand{\Subcase}[2]{\smallskip\paragraph{\textit{\boldmath Subcase {#1}: {#2}.}}\hfil\break\ignorespaces}

\title[Perfect powers in value sets of polynomials]
{Perfect powers in value sets and orbits of polynomials}

\author[A. Ostafe] {Alina Ostafe}
\address{School of Mathematics and Statistics, University of New South Wales, Sydney NSW 2052, Australia}
\email{alina.ostafe@unsw.edu.au}

\author[L. Pottmeyer] {Lukas Pottmeyer}
\address{Fakult\"at f\"ur Mathematik, Universit\"at Duisburg-Essen, 45117 Essen, Germany}
\email{lukas.pottmeyer@uni-due.de}

\author[I. E. Shparlinski] {Igor E. Shparlinski}
\address{School of Mathematics and Statistics, University of New South Wales, Sydney NSW 2052, Australia}
\email{igor.shparlinski@unsw.edu.au}



\begin{document}

\begin{abstract} 
We show the finiteness of perfect powers in orbits of polynomial 
dynamical systems over an algebraic number field.   We also obtain similar results for perfect powers represented by 
ratios of consecutive elements in orbits. Assuming the $abc$-Conjecture for number fields, we obtain a finiteness result for powers in ratios of arbitrary elements in orbits. 
\end{abstract}

\maketitle


\section{Introduction and statements of main results}

\subsection{Motivation}  
Cahn, Jones  and Spear~\cite{CJS} have recently obtained a series of results about the structure 
of intersections of orbits 
of a rational function $\psi \in \LL(X)$   over a  field $\LL$ of characteristic zero,  
 with the image
$\varphi(\LL)$ of  $\LL$  for 
another rational function $\varphi \in \LL(X)$.  
 In the special case  of  $\varphi(X) = X^m$, with a fixed integer $m \ge 2$, this corresponds to the case
of powers in orbits of rational functions, see~\cite[Corollaries~1.6--1.8]{CJS}.  
In particular,
Cahn, Jones  and Spear~\cite[Corollary~1.8]{CJS} give a very explicit characterisation of   polynomials $f(X) \in \LL[X]$ for which for some $\alpha \in \LL$ the 
intersection of the  orbit of $\alpha$  
with the set of $m$-th powers 
$\LL^m$ is  finite.

Here we consider  this question for polynomials  $f(X) \in \KK[X]$  over a number field $\KK$ and 
extend it in two directions, namely,
\begin{itemize}
\item  we consider the union of all orbits over all   $\alpha \in \KK$,
\item we study its intersection with the set of all nontrivial powers of $S$-integers, 
\end{itemize}
see Section~\ref{subsection:notation} for exact definitions.

In particular, we achieve some progress towards Conjecture~\ref{conj:fnklau}  below, which asserts 
that under mild assumptions on $f$ there are only finitely many orbits over 
$\KK$ for which a  ratio of two different iterates is a perfect power. 
We establish this conjecture assuming a $\KK$-rational  $abc$-Conjecture, 
see Theorem~\ref{thm:W K-cond}. 
Furthermore, for a rather large class of polynomials, we establish Conjecture~\ref{conj:fnklau} unconditionally
in Theorem~\ref{thm:set V 0}. 

In fact, we put this question in a more general context of powers in images of polynomials, that is, 
in $f(\KK)$ and reduce it to a much more studied question about powers in the set $f(R_\Scal)$, where $R_\Scal$ 
is a    ring of $\Scal$-integers of $\KK$, see Section~\ref{subsection:notation} for exact definitions. An application of Northcott's Theorem~\cite{Northcott} then allows us to study powers in orbits.

We note that since the proofs of our results rely on applying Siegel's Theorem~\cite[Theorem~D.8.4]{HinSil}, this in turns 
makes us to always assume that $f$ has at least three  simple roots. We also use results of 
B{\'e}rczes,  Evertse and  Gy\H{o}ry~\cite[Theorems~2.1 and~2.2]{BEG} which require 
$\deg f \ge 2$ and also that $f$ has no multiple roots. Thus to combine these two conditions 
we always assume that $\deg f \ge 3$ and $f$ has only simple roots.

\subsection{Notation and conventions}
\label{subsection:notation}
We set the following notation, which remains fixed for the remainder of
this paper:
\begin{notation}
 \item[\textbullet]  $\KK$ is a number field.
  \item[\textbullet]  $\ZZ_\KK$ is the ring of algebraic integers of $\KK$. 
  \item[\textbullet] $M_{\KK}$ is the set of all places of $\KK$ which we partition into the sets $M_{\KK}^0$ and $M_{\KK}^{\infty}$   of all non-Archimedean and Archimedean places of $\KK$, respectively.
     \item[\textbullet]  $\Scal \subseteq M_{\KK}$ is a finite set of places of $\KK$.
     \item[\textbullet] $R_\Scal$ is the ring of $S$-integers of $\KK$.
      \item[\textbullet] $R_\Scal^*$ is the group of $S$-units of $\KK$.
 \item[\textbullet]  $\ov\KK$ is an algebraic closure of $\KK$.
 \item[\textbullet]  $f(X)\in\KK[X]$ is a polynomial of degree $d\ge2$.
 \item[\textbullet]  For $n\ge0$, we write $f^{(n)}(X)$ for the $n$th iterate of~$f$, that is, 
  \[
    f^{(n)}(X) = \underbrace{f\circ f\circ\cdots\circ f}_{\text{$n$ copies}}(X).
  \]
 \item[\textbullet]  For $\a\in\PP^1(\KK)$, we write $\Orbit_f(\a)$
  for the (forward) orbit of~$\a$,  that is, 
 \[
   \Orbit_f(\a) = \left \{ f^{(n)}(\a):~n\ge 0 \right\}.
 \]
 \item[\textbullet] $\Per(f)$ is  the set of periodic points of $f$ in $\ov\KK$,  that is, 
   the set of points $\a\in\ov\KK$ such that $f^{(n)}(\a)=\a$ for some $n\ge 1$. 
 \item[\textbullet] $\PrePer(f)$ is  the set of preperiodic points of $f$ in $\ov\KK$,  that is, 
   the set of points $\a\in\ov\KK$ such that $\Orbit_f(\a)$ is finite.
 \item[\textbullet] $\Wander_\KK(f)$ is the complement of the set $\PrePer(f)$ in $\KK$,   that is,
         the set $\KK\setminus \PrePer(f)$ of $\KK$-rational wandering points for~$f$.
 \item[\textbullet] $\ZZ_{\ge r}$ denotes the set of integers $n \ge r$, where $r$ is a real number.  
 \item[\textbullet]  $\ZZ_{\ne 0,1}$ denotes the set of integers different from $0$ and $1$.
\end{notation}

It is also convenient to define the function
\[
  \logplus t =\log\max\{t,1\}. 
\]
For every $v\in M_{\KK}$ we denote by $\vert . \vert_v$ the corresponding absolute value on $\KK$, normalised so that the absolute logarithmic Weil height
$h:\ov{\KK}\to[0,\infty)$ is defined by
\begin{equation}
\label{eqn:defht}
  h(\b) = \sum_{v\in M_{\KK}} \logplus \left( |\beta|_v \right).
\end{equation}
Every $v\in M_{\KK}^0$ is induced by a prime ideal $\mathfrak{p}_v \subset \ZZ_{\KK}$. For any $\alpha\in\KK$ and any $v\in M_{\KK}^0$ we set $v(\alpha)=\ord_{\mathfrak{p}_v}(\alpha)$. Hence, the map $v:\KK\twoheadrightarrow\ZZ$ is the normalised valuation on $\KK$ corresponding to the non-Archimedean absolute value $\vert . \vert_v$. See~\cite{BoGu,HinSil,Lang} for further details on absolute values and height functions.

\subsection{Main results}
Below we define three sets of ``exceptional
values''. Our goal is to characterise when these sets may be infinite.

 With~$\KK$, $\Scal$ and~$f$ as defined in
Section~\ref{subsection:notation},  we define 
the set
 \[
 \Ucal(\KK,f,\Scal) =
  \left\{ \alpha \in \KK:~ \exists  (\ell,a)  \in  \ZZ_{\ne 0,1}   \times R_\Scal   \text{ such that}\ f(\a)= a^{\ell}    \right\}, 
 \]
 and show its finiteness under certain conditions on $f$, see also Remark~\ref{rem:Cond} 
 about their necessity. 
 
 Additionally, motivated by obtaining a finiteness result for  ratios of elements in orbits, which are perfect powers, we study the finiteness of the set
\begin{align*}
 \Vcal_{m}(\KK&,f,\Scal) =
  \left\{ \begin{aligned} \alpha \in \KK:~ \exists   (n,\ell,a)  \in \{1,\ldots,m\} \times \ZZ_{\ne 0,1} & \times R_\Scal \\
     \text{ such that}\  & f^{(n)}(\a)=  a^{\ell}  \a
    \end{aligned} \right\},
\end{align*}
 where $m\geq 1$ is a fixed rational integer. 
 
 Our main interest in the set $ \Vcal_{m}(\KK,f,\Scal)$   stems from 
Conjecture~\ref{conj:fnklau}  below, which asserts that there are only finitely many orbits over $\KK$ such that the 
 ratio of different iterates is a perfect power. To study this formally we introduce the set
\[
 \Vcal(\KK,f,\Scal)  =
    \left\{\begin{aligned}  \alpha \in \KK:~ \exists  
      (n,\ell,a) &\in\ZZ_{\ge1} \times \ZZ_{\ne 0,1} \times R_\Scal \\
      &\qquad\qquad\quad \text{ such that}\ f^{(n)}(\a) = a^{\ell} \a\\
    \end{aligned} 
    \right\}.
\]

 \begin{remark}  As we have mentioned, our motivation to investigate the set $\Vcal(\KK,f,\Scal) $ stems from the observation that its finiteness, coupled with Northcott's Theorem~\cite{Northcott}, 
see also~\cite[Lemma~2.3]{BOSS}, 
immediately implies the finiteness of the set of $ \alpha \in \KK $, for which the ratio of two elements in
 $\Orbit_f(\a)$ is in the set of powers. In other words,  it implies the finiteness 
 of the set
\begin{align*}
 \widetilde{\Vcal}(\KK&,f,\Scal)  =
    \left\{\begin{aligned}  \alpha \in \KK:~ \exists  
      (n,k&,\ell,a) \in\ZZ_{\ge1} \times \ZZ_{\ge0}\times \ZZ_{\ne 0,1} \times R_\Scal\\
      &\qquad\quad \text{such that}\ f^{(n+k)}(\a) =  a^{\ell} f^{(k)}(\a) \\
    \end{aligned} 
    \right\}.
 \end{align*}
\end{remark}

If $a\in R_\Scal^*$ in the sets above, then finiteness conditions for these sets have been 
given in~\cite[Theorems~1.2, 1.3 and~1.4]{BOSS}. Our results hold for $S$-integers $a$, 
rather than only for $S$-units.

\begin{theorem}
\label{thm:set U} 
Let $f\in\KK[X]$ be of degree $d\ge 3$, having only simple roots. Then for any finite set of places $\Scal$ of $\KK$, 
 the set $ \Ucal(\KK,f,\Scal)$ is  finite.
\end{theorem}

\begin{remark} \label{rem:Cond} 
As noted above, the condition on the polynomial $f$ in Theorem~\ref{thm:set U}  to have \emph{only} simple roots comes from using~\cite[Theorems~2.1 and~2.2]{BEG}, however it is likely that this condition can be relaxed to just imposing that $f$ has at least three distinct simple roots. We note that this latter condition is indeed needed as the following example shows. Let $g\in \ZZ[X]$ be arbitrary, and set $f(X)=2g(X)^2(X^2-1)$. Then, for all of the infinitely many solutions $(r,s)$ of the Pell equation $X^2 -2Y^2 =1$, we have $f(r)= (2g(r) s)^2$. This contradicts the finiteness of $\Ucal(\KK,f,\Scal)$ for any number field $\KK$ and any set of primes $\Scal$. 
Thus we indeed need  $f$ to have at least three simple roots.
\end{remark}

\begin{remark}
It is important to emphasise the difference between the results of  B{\'e}rczes,   Evertse and  Gy\H{o}ry~\cite[Theorems~2.1 and~2.2]{BEG}, 
which apply to  polynomial values with arguments $\alpha  \in R_\Scal$, and Theorem~\ref{thm:set U}, which is based on some additional 
arguments stemming from~\cite{OSSZ1}  and also  on Siegel's Theorem~\cite[Theorem~D.8.4]{HinSil}, 
  and extends this to the finiteness 
of $\alpha \in \KK$, and also allows negative exponents $\ell$. 
\end{remark}

By Northcott's Theorem~\cite{Northcott},  for any  $\beta \in \KK$ there are only finitely many $\alpha \in \KK$ such that 
$\beta \in  \Orbit_f(\a)$. 
Hence,  from Theorem~\ref{thm:set U}, we have the following direct consequence about  powers in orbits.

\begin{corollary}
\label{cor:pow iter}
Let $f\in\KK[X]$ be of degree $d\ge 3$, having only simple roots. Then for any finite set of places $\Scal$ of $\KK$,  there are at most finitely many $\alpha\in\KK$ such that $f^{(n)}(\a) \in R_\Scal^\ell$ for some  $ (n,\ell) \in\ZZ_{\ge1} \times \ZZ_{\ne 0,1}$. 
\end{corollary}

\begin{remark}
We note that 
Theorem~\ref{thm:set U} shows finiteness of the set of tuples 
$(n,\ell,\a,a) \in \ZZ_{\ge 1}\times \ZZ_{\ne 0,1} \times \Wander_\KK(f) \times R_\Scal$ such that 
$ f^{(n)}(\a)= a^{\ell}$. Indeed, by Theorem~\ref{thm:set U} there are finitely many possible values for 
$f^{(n-1)}(\a)$ such that $ f^{(n)}(\a)= a^{\ell}$, and thus, by Northcott's Theorem, finitely many $(n, \a) \in \ZZ_{\ge 1}\times \Wander_\KK(f)$. \end{remark}

We also make: 

\begin{conj}
\label{conj:fnklau}
Let $f\in\KK[X]$ be of degree $d\ge 3$, having only simple roots and such that $0 \not \in \Per(f)$. 
Then for any finite set of places $\Scal$ of $\KK$,  the set $ \Vcal(\KK,f,\Scal)$ is finite.
\end{conj}

We now provide several results towards Conjecture~\ref{conj:fnklau}.
First we consider the set 
$ \Vcal_{m}(\KK,f,\Scal)$, 
which corresponds to the choices $n\leq m$  in the definition 
of $ \Vcal(\KK,f,\Scal)$. 

\begin{theorem}
\label{thm:set V}
Let $m \in\ZZ_{\geq 1}$ and $f\in\KK[X]$ be of degree $d\ge 3$, having only simple roots and such that $f^{(k)}(0)\ne 0$ for all $k\in\{1,\ldots,m\}$. Then for any finite set of places $\Scal$ of $\KK$,  the set $ \Vcal_{m}(\KK,f,\Scal)$ is finite.
\end{theorem}

As in the above, combining Theorem~\ref{thm:set V} with Northcott's Theorem~\cite{Northcott}, we have the following direct consequence about the ratio of two consecutive elements in orbits.

\begin{corollary} 
Let $f\in\KK[X]$ be of degree $d\ge 3$, having only simple roots and such that $f(0)\ne 0$. Then for any finite set of places $\Scal$ of $\KK$, there are at most finitely many $\alpha\in\KK$ such that $f^{(n+1)}(\a)/f^{(n)}(\a)\in R_\Scal^\ell$  for some  $ (n,\ell) \in\ZZ_{\ge1} \times \ZZ_{\ne 0,1}$.
\end{corollary}

\begin{remark} 
\label{rem: V vs W} The proof of  Theorem~\ref{thm:set V}  splits into several cases, depending on some 
additional assumptions on $\alpha$ and $a$.  In all but one of 
these cases (where we have to assume that $n$ is bounded), 
we actually give a finiteness result for 
their contribution to the set $\Vcal(\KK,f,\Scal)$. 
This exceptional case is a reason why we have a complete proof of finiteness only of the set $ \Vcal_{m}(\KK,f,\Scal)$, rather than resolve Conjecture~\ref{conj:fnklau} in full.
\end{remark}

We now produce an infinite class of polynomials for which Conjecture~\ref{conj:fnklau} holds. 

\begin{theorem}
\label{thm:set V 0}
Let $f\in\KK[X]$ be of degree $d\ge 3$, having only simple roots and such that $0 \in  \PrePer(f) \setminus \Per(f)$. 
Then for any finite set of places $\Scal$ of $\KK$, the set $ \Vcal(\KK,f,\Scal)$ is finite. 
\end{theorem}

For example,  the classes of polynomials 
\[
f(X) = X^n(X^m-1)+\zeta \qquad\mbox{and}\qquad f(X) = X^k(X-b)+b
\]
satisfy the conditions of  Theorem~\ref{thm:set V 0}, for all $n,m\ge 1$ with $n+m\ge 3$ and a root of unity $\zeta$ of order dividing $m$, and
 $k \ge 2$, $b\in \KK\setminus \{0\}$ with $b^k k^k\neq (k+1)^{k+1}$, respectively.
Indeed, it is immediate to check that for these polynomials we have $0\in \PrePer(f) \setminus \Per(f)$. Applying the well-known formula for the discriminant of a trinomial 
(cf.~\cite[Theorem~2]{Swan}) yields that the discriminant of these polynomials is not zero. Hence, all these polynomials have only simple roots.

Next, we prove that Conjecture~\ref{conj:fnklau} follows from the $abc$-Conjecture for the number 
field $\KK$, see~\cite[Chapter~14]{BoGu}.  
To formulate the $abc$-Conjecture for $\KK$, for each $v \in M_{\KK}^0$ we fix a uniformizer $\pi_v \in \KK$, which is just an element in $\KK$ with $v(\pi_v)=1$.

\begin{conj}[$\KK$-rational $abc$-Conjecture]\label{Kabc}
For every $\varepsilon>0$, there exists a constant $C(\varepsilon)$ such that for all $\alpha\in\KK\setminus\{0,1\}$ we have
\begin{equation}\label{eq:abc}
\begin{split} 
(1-\varepsilon) h(\alpha) \leq \sum_{\substack{v(\alpha)>0 \\ v\in M_{\KK}^0}} \log{\left| \frac{1}{\pi_v}\right|_v} &+ \sum_{\substack{v(1-\alpha)>0 \\ v\in M_{\KK}^0}} \log{\left| \frac{1}{\pi_v}\right|_v}\\
& + \sum_{\substack{v(1/\alpha)>0 \\ v\in M_{\KK}^0}} \log{\left| \frac{1}{\pi_v}\right|_v} + C(\varepsilon).
\end{split} 
\end{equation}
\end{conj}  

\begin{remark}
In the case $\KK=\QQ$, the validity of Conjecture~\ref{Kabc}  for any   $\varepsilon>0$ is equivalent to the classical $abc$-Conjecture of Masser and Oesterl\'e. Namely, for any $\varepsilon>0$, there exists a constant $C(\varepsilon)$, such that for all pairwise coprime $a,b,c \in \ZZ_{\geq 1}$, with $a+b=c$, we have 
\[
\prod_{\substack{p \mid abc\\p~\text{prime}}} p \ge C(\varepsilon) c^{1-\varepsilon}.
\]
For   further information on Conjecture~\ref{Kabc} we refer to~\cite[Chapter~14]{BoGu}. 
\end{remark}

\begin{theorem}
\label{thm:W K-cond} 
Let $f\in \KK[X]$ be of degree $d\ge 3$, having only simple roots and such that $0 \not \in \Per(f)$. 
Assuming the validity of the $\KK$-rational $abc$-Conjecture~\ref{Kabc}, for any finite set of places $\Scal$ of $\KK$, the set $\Vcal(\KK,f,\Scal) $ is finite. 
\end{theorem}

In other words,  Theorem~\ref{thm:W K-cond}  can be informally expressed as 
$$
\text{Conjecture~\ref{Kabc}}  \quad \Longrightarrow \quad \text{Conjecture~\ref{conj:fnklau}}.
$$

The assumption $0\not\in \Per(f)$ is needed in the proof of Theorem~\ref{thm:W K-cond}. However, we do not have a counterexample for Conjecture~\ref{conj:fnklau} if we only assume $f(0)\neq 0$. The latter assumption is indeed necessary, since the same argument as in Remark~\ref{rem:Cond} shows that for $f(x)=2x(x^2 -1)$ the set $\Vcal(\QQ,f,\Scal)$ is infinite.

\begin{remark}
Examining the proof of Theorem~\ref{thm:W K-cond} one can see that for a fixed polynomial $f\in \KK[x]$ satisfying the assumptions from Theorem~\ref{thm:W K-cond}, we do not need the full power  of the $\KK$-rational $abc$-Conjecture to prove the finiteness of $\Vcal(\KK,f,\Scal)$. We  prove that for any such $f$, there exists an $\varepsilon(f)>0$ such that if~\eqref{eq:abc} holds for $\varepsilon(f)$, then $\Vcal(\KK,f,\Scal)$ is finite. See also Remark~\ref{rem:eps Bely} after the proof of Theorem~\ref{thm:W K-cond} for more details.
\end{remark}

\subsection{Ideas behind the proofs} 
To simplify the exposition we outline the main ingredients of our arguments in the case of $\KK = \QQ$, and thus instead of the language of valuations, we simply talk about prime divisors. 
One of the most powerful tools we use is a result of B{\'e}rczes, Evertse 
and  Gy\H{o}ry~\cite{BEG}, which gives 
a very explicit form of a result of Schinzel and Tijdeman~\cite{SchTij} on the finiteness of  perfect powers among the integral values of polynomials over $\ZZ$ modulo a finite set of primes. Say, in Theorem~\ref{thm:set U} this means, using~\cite{BEG}, that for $\ell\ge 2$ and any  $\alpha\in \QQ$ with a denominator composed out of a fixed set of primes  the power $\ell$ is uniformly bounded  which in turn easily implies the boundedness of the set of possible
solutions $ (\a, \ell, a)$.  If $\ell<0$, then we show that the numerator of $f(\alpha)$ is  composed out of a fixed set of primes and the result is a simple consequence of the celebrated result of Siegel~\cite[Theorem~D.8.4]{HinSil}.

For the proof of Theorem~\ref{thm:set V 0}, we note that $f^{(n)}(\a) =  a^{\ell} \a$ means that $\a$ ``almost'' divides $f^{(n)}(0)$. Now, if $f^{(n)}(0)\ne 0$, but takes just finitely many different values, also $\alpha $ can take just finitely  many different values.

In the proof of Theorem~\ref{thm:W K-cond} we first handle small values of $n$. Here we proceed as we described before and conclude that, since $\a$ ``almost'' divides $f^{(n)}(0)$, we are able to apply~\cite{BEG}. When $n$ is large we need the power of the $abc$-Conjecture to conclude, as in~\cite{Gran}, that $\alpha$ must be bounded as polynomial values can generally not be multiplicatively close to perfect powers. 

Since Siegel's Theorem enters our argument,  Theorems~\ref{thm:set U}, \ref{thm:set V}, and \ref{thm:W K-cond}
  are not effective.

\begin{remark}
It is easy to see that our finiteness results are also true for shifts by $\Scal$-units, that is, we can replace the set $\Ucal(\KK,f,\Scal)$ by the set
\begin{equation}\label{eq:U and u}
\left\{ \alpha \in \KK:~ \exists  (u,\ell,a)  \in  R_{\Scal}^* \times \ZZ_{\ne 0,1}   \times R_\Scal   \text{ such that}\ f(\a)= u a^{\ell}    \right\},
\end{equation}
and similarly for $\Vcal_m(\KK,f,\Scal)$ and $\Vcal(\KK,f,\Scal)$. Indeed, if $\alpha\in R_\Scal$, then by~\cite[Lemma~2.8]{BOSS} we may assume that $\ell$ must be bounded from above. Hence, after replacing $\KK$ by a larger number field, we may assume that all $\Scal$-units are $\ell$-th powers for all possible (finitely many) values of $\ell$. Hence, extending if necessary the ground field $\KK$ to contain 
all necessary $\ell$-roots of all generators of $R_{\Scal}^*$,  we see that 
the finiteness of the subset of the set~\eqref{eq:U and u} with  $\alpha\in R_\Scal$ follows directly from  Theorems~\ref{thm:set U}.  We also note that in  the case $\alpha\not\in R_\Scal$ the proof of Theorem~\ref{thm:set U} holds without any change for the equation $f(\a)= u a^{\ell}$ in~\eqref{eq:U and u}  thus concluding the proof of the finiteness of the set~\eqref{eq:U and u}.

Furthermore, trivial changes in the proofs of   Theorems~\ref{thm:set V}, \ref{thm:set V 0} and~\ref{thm:W K-cond} below give the finiteness of the similar  generalisations of  $\Vcal_m(\KK,f,\Scal)$ and $\Vcal(\KK,f,\Scal)$.
\end{remark}

\section{Proofs of main results} 

\subsection{Preliminary discussion} 

As usual, we say that
a polynomial
\[
f(X)=c_0+c_1X+\cdots+c_dX^d
\]
has bad reduction at~$v\in M_{\KK}^0$
if either $v(c_i)<0$ for some~$i$ or if~$v(c_d)>0$; otherwise
we say it has good reduction.  We   let 
\[
  \Scal_{f} =  \Scal \cup M_\KK^\infty\cup
   \{v\in M_\KK^0 : \text{$f$ has bad reduction at $v$} \}.
\]
In particular,  for all $v\not\in\Scal_{f}$ we have 
\begin{equation}
\label{eq:qinSf}
|c_d|_{v}=1 \qquad\mbox{and}\qquad |c_i|_{v}\le 1, \qquad i=0,\ldots,d-1.
\end{equation}
 
It is easy to see  that if~$f$ has
good reduction at $v$, then so do all of its iterates; in fact this  is 
also true even for rational functions,
see~\cite[Proposition~2.18(b)]{Silv07}. 
Hence
\begin{equation}
\label{eq:SmS}
  \Scal_{f^{(m)}} \subseteq \Scal_{f}, \quad \text{for all $m\ge1$.}
\end{equation}

We let
\[
  \RSf = \left\{ \vartheta \in \KK:~\text{$ v(\vartheta)\ge 0$ for all $v \not\in \Scal_{f}$} \right\}
\]
be the ring of $\Scal_{f}$-integers in $\KK$, and $ \RSf^*$
denotes the group of $\Scal_{f}$-units in $\KK$. Clearly $R_\Scal\subseteq \RSf$.

\subsection{Proof of Theorem~$\ref{thm:set U}$} 
We replace the set $R_\Scal$ with  $ \RSf$ and thus   investigate the equation 
\begin{equation}
\label{eq:ful}
f(\alpha) = a^\ell, \qquad (\a, \ell, a) \in  \KK \times\ZZ_{\ne 0,1} \times  \RSf .
\end{equation} 
If $\a \not\in \RSf$, then there is some non-Archimedean valuation $v\in M_{\KK}\setminus \Scal_{f}$ such that $v(\a)<0$. Hence, recalling~\eqref{eq:qinSf}, we see from the basic properties of non-Archimedean valuations  that $v(f(\alpha)) = d  v(\alpha)$ (see also  the proof of~\cite[Theorem~4.11]{OSSZ1}). Now, if $(\a,\ell,a)$ is a solution to \eqref{eq:ful} with $\a\not\in \RSf$ then
\[
0 > dv(\a)=v(f(\a))=v(a^{\ell})=\ell v(a).
\]
Since $a\in \RSf$, it follows that $\ell < 0$. We now distinguish between two cases.

\Case{A}{$\ell \geq 2$} By our preliminary discussion, we know that in this case we have $\alpha \in \RSf$. Hence, this case follows directly from~\cite{BEG}. Indeed, by~\cite[Theorem~2.3]{BEG} (see also~\cite[Lemma~2.8]{BOSS}) the exponent $\ell\ge 2$ is bounded  by a constant depending only on $\KK$, $f$ and $\Scal$. Thus, we can consider $\ell$ to be fixed, and since $\deg f\ge 3$, we can apply~\cite[Theorems~2.1 and~2.2]{BEG}   to conclude that $h(\a)$ is bounded by a constant depending only on $\KK$, $f$ and  $\Scal$. Now, since $\a\in\KK$, then Northcott's Theorem~\cite[Theorem~3.7]{Silv07} tells us that there are finitely many such $\a$.

\Case{B}{$\ell< 0$} Let us define the rational function $g(X)=f(X)^{-1}$. Then, since $f$ has at least three simple roots, the function $g$ has at least three distinct poles. 

We have a solution $f(\a)=a^{\ell}$ with $a\ne 0$, if and only if 
\[
g(\a)=f(\a)^{-1}=a^{-\ell}.
\]
Since $-\ell>0$, $a^{-\ell}\in\RSf$, and thus we can apply Siegel's Theorem~\cite[Theorem~D.8.4]{HinSil}   
to conclude that there are finitely many $\alpha\in\KK$ such that $g(\alpha)\in\RSf$. This concludes also this case.

\subsection{Proof of Theorem~\ref{thm:set V}}
As we have explained in Remark~\ref{rem: V vs W},  
we  follow the proof in the general case of proving Conjecture~\ref{conj:fnklau}, except for one case which breaks down, and thus prove this case only for  $n\leq m$, for some fixed $m\in\ZZ_{\geq1}$, which  concludes only the proof of Theorem~\ref{thm:set V}.
In particular, we consider the equation
\begin{equation}
\label{eq: fna}
f^{(n)}(\alpha)=a^{\ell}\alpha.
\end{equation} 

Moreover, since by Northcott's Theorem~\cite{Northcott} (see also~\cite[Theorem~3.12]{Silv07}), the set $\PrePer(f)\cap\KK$ is finite, we need only to prove finiteness of the set $ \Vcal(\KK,f,\Scal)\cap \Wander_\KK(f)$.

We split now the proof into two cases depending on $\ell$ being positive and negative, and then some further subcases.

\Case{A}{$\ell\ge 2$}  
Assume there is a solution $\a\in \Wander_\KK(f)$ to~\eqref{eq: fna}.  Since $\ell\ge 0$,  for all $v\in M_{\KK} \setminus \Scal_{f}$, one has
\begin{equation}
\label{eq:f n v}
|f^{(n)}(\alpha)|_{v}\le |\alpha|_{v}.
\end{equation}
If there exists $v\in M_{\KK} \setminus \Scal_{f}$ such that $|\alpha|_v>1$, then~\eqref{eq:qinSf} and a simple valuation computation show that
$$
|f^{(n)}(\alpha)|_v=|\alpha|_v^{d^n}.
$$
This together with~\eqref{eq:f n v} and the fact that $d\ge 2$, leads to a contradiction. Thus, for any $v\in M_{\KK} \setminus \Scal_{f}$ we have $|\alpha|_v\leq 1$, which implies that 
\[
\alpha\in \RSf.
\]
We now consider the following three subcases. 

\Subcase{A.1}{$a\in \RSf^*$}
The finiteness in this case follows directly from~\cite[Theorem~1.3]{BOSS}.

\Subcase{A.2}{$a\in\RSf\setminus \RSf^*$ and $\alpha\in \RSf^*$} 
Precisely as in \emph{Case A}  in the proof of Theorem~\ref{thm:set U} (but using~\cite[Lemma~2.8]{BOSS} instead of~\cite[Theorem~2.3]{BEG}) we may assume that $\ell$ is fixed. Since $\RSf^*$ is finitely generated,  we may take a finite field extension $\LL/\KK$ such that $\LL$ includes the $\ell$-th roots of all the generators of $\RSf^*$. We also extend the set $\Scal_f$ to $\widetilde{\Scal}_f$ including all the places of $\LL$ staying over the places in $\Scal_f$. Now, finiteness follows immediately from Corollary~\ref{cor:pow iter}.

\Subcase{A.3}{$a\in\RSf\setminus \RSf^*$ and $\alpha\in \RSf\setminus\RSf^*$}
In this case, since $\alpha\in \RSf\setminus\RSf^*$, there exists $v\not\in\Scal_{f}$ such that $v(\alpha)>0$. Since $a\in\RSf$, we also have $v(a)\ge 0$, and thus $v(f^{(n)}(\alpha))>0$ (since $\ell>0$).

Now, let us write 
\[
f^{(n)}(X)=\sum_{i=0}^{d^n} c_{i,n} X ^i,
\]
with $v( c_{i,n})\ge 0$, which follows 
from~\eqref{eq:qinSf} and~\eqref{eq:SmS}. 
Thus $v(c_{i,n}\alpha^i)>0$ for all $i \in \{1,  \ldots, d^n\}$ and 
\begin{equation}
\label{eq:vq_f0}
\begin{split}
v(f^{(n)}(0))&=v(c_{0,n})=v\left(f^{(n)}(\alpha)-\sum_{i=1}^{d^n}c_{i,n}\alpha^i\right)\\
&\ge \min \{v(f^{(n)}(\alpha)),\min_{i=1,\ldots,d^n}v(c_{i,n}\alpha^i)\}>0.
\end{split}
\end{equation}
Thus, for any $v\not\in\Scal_{f}$ such that $v(\alpha)>0$, one has $v(f^{(n)}(0))>0$. 

This is where we do not know how to conclude the proof in full generality and thus {\it for the rest of Subcase~A.3 only\/} we assume  
\[
n\leq m.
\]
By~\eqref{eq:vq_f0}, for any $v\not\in\Scal_{f}$ such that $v(\alpha)>0$, one has $v(f^{(n)}(0))>0$. However, since $f^{(n)}(0)\neq 0$ for all $n\leq m$, there are at most finitely many $v\in M_{\KK}^0$ such that $v(f^{(n)}(0))>0$ for some $n\in\{1,\ldots,m\}$.
Thus, extending $\Scal_{f}$ to include all these places  and denoting this new set by $\Tcal_{f,m}$, we can conclude that $\alpha\in R_{\Tcal_{f,m}}^*$. The finiteness conclusion follows now  as in {\it Subcase~A.2\/} applied with $\Tcal_{f,m}$ instead of $\Scal_{f}$ (and noting  that $\RSf\subseteq R_{\Tcal_{f,m}}$). 
This concludes thus the finiteness of the set $ \Vcal_{m}(\KK,f,\Scal)$.

\Case{B}{$\ell<0$} 
We continue now the proof for arbitrary $n\ge 1$.

Since $f$ is a polynomial of degree $d\ge 3$ with only simple roots, $f$ is not of the form $cX^d$ and moreover $0$ is not an exceptional point for $f$ (if $0$ would be an exceptional point, then the cardinality of the backward orbit of $0$ would be $1$ or $2$, see for example~\cite[Theorem~1.6]{Silv07}, which is impossible).

We study the finiteness of the set of  elements   $\a\in  \Wander_\KK(f)$ such that
\begin{equation}
  \label{eqn:fnkvafnavvnotS}
      |f^{(n)}(\a)|_v = |a|_v^{\ell} |\a|_v,  \qquad \forall\, v\in M_\KK \setminus \Scal_{f},
\end{equation}
 for some $(n,\ell,a) \in\ZZ_{\ge1}\times \ZZ_{<0} \times (\RSf\setminus\{0\})$.

We now proceed as in the proof of~\cite[Theorem~1.3]{BOSS} and we  indicate only what is new. For an arbitrary choice of~$\varepsilon$, to be specified later, we
let $C_3(\KK,\Scal_{f},f,\varepsilon)$ be the constant from~\cite[Lemma~2.5]{BOSS}, and we split the proof into two
cases, depending whether  $n$ is large or small. 

\Subcase{B.1}{$n\ge C_3(\KK,\Scal_{f},f,1/3)$}
In this case, by~\cite[Lemma~2.5]{BOSS} applied with  $\varepsilon = 1/3$, we see that $(n,\alpha)$
satisfies
\begin{equation}
  \label{eqn:lpfnkleehfnk}
  \sum_{v\in \Scal_{f}}
  \logplus\left(|f^{(n)}(\a)|_v^{-1}\right) \le  \frac{1}{3} \hhat_f \left(f^{(n)}(\a)\right),
\end{equation}
where $\hhat_f$ is the canonical height associated to $f$, see~\cite[Section~3.4]{Silv07} 
for a definition and  standard properties. 

Since $h(\g)=h(\g^{-1})$ and using~\eqref{eqn:defht}, we compute
\begin{align*}
  h\left(f^{(n)}(\a)\right)
 &= h\left(f^{(n)}(\a)^{-1}\right)  
 = \sum_{v\in M_\KK} \logplus\left(|f^{(n)}(\a)|_v^{-1}\right)\\
   &= \sum_{v\in \Scal_{f}} \logplus\left(|f^{(n)}(\a)|_v^{-1}\right)  +  \sum_{v\in M_\KK \setminus \Scal_{f}} \logplus\left(|f^{(n)}(\a)|_v^{-1}\right).
\end{align*}
Now, using~\eqref{eqn:fnkvafnavvnotS} and~\eqref{eqn:lpfnkleehfnk} and the fact that $\ell<0$ (and thus $|a^{-\ell}|_v\le 1$ for all $v\in M_\KK \setminus \Scal_{f}$), we have
\begin{align*}
  h\left(f^{(n)}(\a)\right) 
   &\le  \frac{1}{3}  \hhat_f \left(f^{(n)}(\a)\right) +  \sum_{v\in M_\KK \setminus \Scal_{f}}\logplus\left(|a^{\ell}\a|_v^{-1}\right)
 \\
 &\le  \frac{1}{3}  \hhat_f \left(f^{(n)}(\a)\right) +    h\left(\a^{-1}\right)  =  \frac{1}{3}  \hhat_f \left(f^{(n)}(\a)\right) +  h\left(\a\right).
\end{align*}

From now on the proof goes word by word as in  the proof of~\cite[Theorem~1.3]{BOSS} with  $\rho=1$ and $k=0$ (where also  a somewhat arbitrary 
value $\varepsilon=1/3$ has been used). 
This implies also the finiteness of the set $ \Vcal(\KK,f,\Scal)$ in this case.

\Subcase{B.2}{$n< C_3(\KK,\Scal_{f},f,1/3)$} 
Let $g(X)=f^{(n)}(X)/X$, and we note that $g$ has at least three nonzero distinct roots, which follows immediately from the fact that $f$ has this property.   Since $n$ is bounded, proving finiteness of the set $\Vcal(\KK,f,\Scal)$ in this case reduces to proving finiteness of $\a\in\KK$ such that 
$g(\alpha)=a^{\ell}$ for some $(\ell,a) \in  \ZZ_{\ne 0,1} \times R_\Scal$.
This follows exactly as in the proof of Theorem~\ref{thm:set U}, \emph{Case B}, applying Siegel's Theorem~\cite[Theorem~D.8.4]{HinSil}.  
Indeed, we follow the proof of \emph{Case~B} of Theorem~\ref{thm:set U} above with $f(X)$ replaced by 
the rational function $g(X)$, and apply Siegel's Theorem to the function $G(X)=g(X)^{-1}$ (taking 
also into account that $f^{(n)}(0)\ne 0$) to conclude that $G(\KK)\cap R_{\Scal_f}$ is finite.

\subsection{Proof of Theorem~\ref{thm:set V 0}}
We recall that in the proof of Theorem~\ref{thm:set V} only {\it Subcase~A.3\/} requires the 
assumption that $n$ is bounded. Hence we consider only this case. Recall that this assumption appears after it has been 
shown that  for any $v\not\in\Scal_{f}$ such that $v(\alpha)>0$, one has $v(f^{(n)}(0))>0$. Since  
$0 \in  \PrePer(f) \setminus \Per(f)$ we see that there are only finitely many $v\in M_{\KK}^0$ with this property. 
Hence $\alpha \in \RTf^*$ for some finite set $\Tcal_{f}$ depending only on $\Scal$ and $f$.  
We now proceed as in {\it Subcase~A.2\/}  
in the proof of Theorem~\ref{thm:set V}  and obtain the desired result.

 \subsection{Proof of Theorem~\ref{thm:W K-cond}}

All constants in this proof may depend on the fixed number field $\KK$, even when we do not explicitly state this dependence.

After enlarging the set $\Scal$ to $\Scal_{f}$, we may assume that $f\in R_{\Scal_{f}}[x]$. 

Again, we only have to consider  {\it Subcase~A.3\/} from the proof of Theorem~\ref{thm:set V}. Hence, we have to prove the finiteness of
\begin{equation}\label{finite set W}
\left\{\begin{aligned}  \alpha \in R_{\Scal_{f}}:~ \exists  
      (n,\ell,a) &\in\ZZ_{\geq m} \times \ZZ_{\geq 2} \times R_{\Scal_{f}} \\ &  \qquad \qquad \text{ such that}\ f^{(n)}(\a) = a^{\ell} \a\\
    \end{aligned} 
    \right\},
\end{equation}
for some fixed positive integer $m$. 

It follows from~\cite[Theorem~3.11]{Silv07} that there exists a constant $C_1(f)$ such that for all $k\geq 1$ and all $\alpha\in\KK$ we have
\begin{equation}\label{eq:height equiv}
d^k h(\alpha)- d^k C_1(f) \leq h(f^{(k)}(\alpha)) \leq d^k h(\alpha) + d^k C_1(f).
\end{equation}
In particular, the set of $\alpha\in\KK$ such that $f^{(k)}(\alpha)=0$ for some $k$ is a set of bounded height, and hence it is finite. Therefore, we may assume without loss of generality that $f^{(k)}(\alpha) \neq 0$ for all $k\geq 1$.

Let us fix some further notations in order to apply the $\KK$-rational $abc$-Conjecture~\ref{Kabc}. If $D$ is a divisor on $\PP^1_{\KK}$, then $h_D$ denotes the height associated to $D$. This is, if  $(\lambda_{D,v})_{v\in M_{\KK}}$ is a collection of local heights associated to $D$, then 
\[
h_D(P)=\sum_{v\in M_{\KK}} \lambda_{D,v} (P)
\] 
for all $P\in\PP^1(\overline{\KK})\setminus\Support(D)$, where  $\Support(D)$ 
is the support of $D$. 

A height associated to the canonical divisor of $\PP^1_{\KK}$ can be chosen to be $-2h$, where 
\[
h([x_0:x_1])=\sum_{v\in M_{\KK}}\max\{\left| x_0\right|_v, \left| x_1\right|_v\}
\] 
is the standard height on $\PP^1(\overline{\KK})$ (for example, see~\cite[Example~14.4.4]{BoGu}). 
Therefore,  we have the following link between Conjecture~\ref{Kabc} and Vojta's conjectured height inequality
 (see~\cite[Theorem~14.4.16  and Remark~14.4.17]{BoGu}):

Let $D$ be a reduced divisor on $\PP^1_{\KK}$. We see that if the $\KK$-rational $abc$-Conjecture~\ref{Kabc} is true, then for all $\varepsilon>0$ there exists a constant $C(\varepsilon,D)$ only depending on $\varepsilon$ and $D$ such that for all $P\in\PP^1(\overline{\KK})\setminus\Support(D)$ one has
\begin{equation}\label{eq:Vojta}
h_D(P) -(2+\varepsilon)h(P)\leq \sum_{\substack{\lambda_{D,v}(P) >0 \\ v\in M_{\KK}^0}}\log\left| \frac{1}{\pi_v}\right|_v +C(\varepsilon,D).
\end{equation}
Note, that the contribution of the places from the finite set $\Scal_f$ to the sum on the right hand side is bonded by a constant. Hence, applying \eqref{eq:Vojta} to $D=\Div(F)$, where $F(X,Y)=Y^{d+1}f(X/Y)\in\RSf[X,Y]$, yields
\begin{equation}\label{eq:Granville}
(d-1-\varepsilon)h(\alpha)\leq \sum_{\substack{v(f(\alpha))>0 \\ v\in M_{\KK}^0\setminus \Scal_{f}}}\log\left| \frac{1}{\pi_v}\right|_v +C_2(\varepsilon,f,\Scal_{f}), \quad \forall~ \alpha \in  R_{\Scal_{f}},
\end{equation}
where $C_2(\varepsilon,f,\Scal_{f})$ is a constant only depending on $\varepsilon$, $f$ and $\Scal_f$. That the $abc$-conjecture implies~\eqref{eq:Granville} is well known. In particular, in the case $\KK=\QQ$ this statement has been used by Granville~\cite{Gran} to count squarefree values of integer polynomials. The implication for number fields is due to Elkies~\cite[Equation~(26)]{Elkies}. In the Elkies inequality~\cite[Equation~(26)]{Elkies}, the term  $(d-1-\varepsilon)$ in~\eqref{eq:Granville} is replaced by $(d-2-\varepsilon)$. This bound comes from using the usual homogenisation $Y^{d}f(X/Y)$ of $f$ instead of $F(X,Y)=Y^{d+1}f(X/Y)$. The same homogenisation has also been used by Granville~\cite{Gran}.

Now, we calculate an upper bound for 
\[
\sum_{\substack{v(f^{(n)}(\alpha))>0 \\ v\in M_{\KK}^0\setminus \Scal_{f}}}\log\left| \frac{1}{\pi_v}\right|_v = \sum_{\substack{\left| f^{(n)}(\alpha)\right|_v<1 \\ v\in M_{\KK}^0\setminus \Scal_{f}}}\log\left| \frac{1}{\pi_v}\right|_v,
\] 
if $\alpha\in \RSf$ satisfies $f^{(n)}(\alpha)=a^{\ell}\alpha\neq 0$ for some $n,\ell\geq 2$ and $a\in R_{\Scal_{f}}$.  In this case it is
\begin{align*} 
\label{eq:conductor vs height}
  \sum_{\substack{\left| f^{(n)}(\alpha)\right|_v <1  \\ v\in M_{\KK}^0 \setminus \Scal_{f}}} \log\left| \frac{1}{\pi_v}\right|_v
& =\sum_{\substack{\left| a^{\ell}\alpha \right|_v <1  \\ v\in M_{\KK}^0 \setminus \Scal_{f}}}\log\left| \frac{1}{\pi_v}\right|_v 
=   \sum_{\substack{\left| a\alpha \right|_v <1  \\ v\in M_{\KK}^0 \setminus \Scal_{f}}}\log\left| \frac{1}{\pi_v}\right|_v \\
& \leq  \sum_{v\in M_{\KK}^0 \setminus \Scal_{f}}\logplus \left| (a\alpha)^{-1}\right|_v  
=   \frac{1}{\ell} \sum_{v\in M_{\KK}^0 \setminus \Scal_{f}}\logplus \left| (a\alpha)^{-\ell}\right|_v.
\end{align*}
Hence we now arrive to the inequality
\begin{equation}
\begin{split} 
\label{eq:conductor vs height}
  \sum_{\substack{\left| f^{(n)}(\alpha)\right|_v <1  \\ v\in M_{\KK}^0 \setminus \Scal_{f}}}  \log\left| \frac{1}{\pi_v}\right|_v
&  \leq   \frac{1}{\ell} h((a\alpha)^{-\ell}) =\frac{1}{\ell}h((a\alpha)^{\ell}) \\
& =\frac{1}{\ell} h(a^{\ell}\a \cdot \a^{\ell-1})  
\leq    \frac{1}{\ell} h(f^{(n)}(\a)) + \frac{1}{\ell}h(\a^{\ell-1})\\
&  \leq  \frac{1}{2} h(f^{(n)}(\a)) + h(\a).
\end{split}
\end{equation}

Now, applying~\eqref{eq:height equiv} to the inequality \eqref{eq:conductor vs height}, we obtain
\begin{equation}\label{eq:conductor upper bound}
\sum_{\substack{\left| f^{(n)}(\alpha)\right|_v <1  \\ v\in M_{\KK}^0 \setminus \Scal_{f}}}  \log\left| \frac{1}{\pi_v}\right|_v \leq \left(\frac{1}{2}d^n +1 \right) h(\alpha)+\frac{1}{2}d^n C_1(f).
\end{equation}

Combining~\eqref{eq:conductor upper bound} with~\eqref{eq:Granville} yields
\[
(d-1-\varepsilon)h(f^{(n-1)}(\alpha)) \leq \left(\frac{1}{2}d^n +1 \right) h(\alpha)+\frac{1}{2}d^n C_1(f) +C_2(\varepsilon,f,\Scal_{f}).
\]
A further application of~\eqref{eq:height equiv} implies that 
\begin{align*}
\left( \frac{1}{2} d^n - (1+\varepsilon)d^{n-1} -1 \right) h(\alpha) - \left( \frac{3}{2} d^n - (1+\varepsilon)d^{n-1} \right)&C_1(f)\\ &\leq C_2(\varepsilon,f,\Scal_{f}).
\end{align*}
It follows that
\[
\left( \frac{1}{2} d^n - (1+\varepsilon)d^{n-1} -1 \right) \cdot \left(h(\alpha) - \frac{3d/2-(1+\varepsilon)}{d/2 -(1+\varepsilon)} C_1(f) \right)
\]
is bounded from above independently on $\alpha$ and $n$. Since we assume that $\varepsilon < d/2 -1 $, the factor $(d^n/2-(1+\varepsilon)d^{n-1}-1)$ tends to infinity as $n$ does. Hence, either $n$ is bounded independently on $\alpha$, or $h(\alpha)$ is bounded independently on $n$. In the first case, finiteness of the set in~\eqref{finite set W} follows immediately from Theorem~\ref{thm:set V}. In the second case the claimed finiteness follows as there are only finitely many points of bounded height and bounded degree.

\begin{remark}
\label{rem:eps Bely}
We conclude with a remark on the precise dependence of Theorem~\ref{thm:W K-cond} on the $abc$-conjecture. The $\KK$-rational $abc$-conjecture implies~\eqref{eq:Vojta} in the following way: If inequality~\eqref{eq:abc} holds for some 
fixed $\varepsilon>0$, then~\eqref{eq:Vojta} holds for 
$\widetilde \varepsilon=n\varepsilon$, where $n$ only depends on the degree of a \emph{Belyi map defined over $\KK$ associated to $\Support(D)$}. Such a Belyi map can always be chosen as a polynomial defined over $\QQ$, and an explicit bound on its degree has been 
calculated in~\cite{Kha}. Then, $n$ is twice the degree of this map. For a proof of this statement, the interested reader may follow the proof 
of~\cite[Theorem~14.4.16: (a) $\Rightarrow$ (b)]{BoGu}. 
We conclude that there exists an (effectively computable) $\varepsilon(f)$, depending solely on $f$, 
such that~\eqref{eq:abc} for $\varepsilon(f)$ implies~\eqref{eq:Granville} for $\widetilde \varepsilon=d/2 -1.0001$. As we have just seen, this again  implies the finiteness of $\Vcal(\KK,f,\Scal_{f})$.
\end{remark}

\section*{Acknowledgement}

The authors are grateful to Attila B\' erczes and Florian Luca for useful discussions and comments on an initial draft of the paper. The authors are very grateful to the anonymous referee for valuable comments and suggestions, which simplified some proofs and improved the presentation of the paper. We also thank Conrad Martin and the referee for noticing a computational error in the initial proof of Theorem~\ref{thm:set U}.

The work of A.O. and I.S.   was  supported in part  by the Australian Research Council Grant~DP180100201.

\end{document}